\documentclass[12pt] {article}

  \usepackage{amsmath}
\usepackage{verbatim}

\newtheorem{example}{Example}[section]
\newtheorem{theorem}{Theorem}[section]
\newtheorem{lemma}{Lemma}[section]
\newtheorem{corollary}{Corollary}[section]
\newtheorem{proposition}{Proposition}[section]
\newtheorem{remark}{Remark}[section]
\newcommand{\eqnsection}{
   \renewcommand{\theequation}{\thesection.\arabic{equation}}
   \makeatletter
   \csname @addtoreset\endcsname{equation}{section}
   \makeatother}
 
\def \ov{\overline}

\def \be{\begin{equation}}
\def \ee{\end{equation}}
\def \bt{\begin{theorem}} 
\def \et{\end{theorem}}
\def \bl{\begin{lemma}} 
\def \el{\end{lemma}}
\def \bea{\begin{eqnarray}}
\def \eea{\end{eqnarray}}
\def \bas{\begin{eqnarray*}}
\def \eas{\end{eqnarray*}}

\def \al{\alpha}
\def \bb{\beta}
\def \ga{\gamma}
\def \Ga{\Gamma}
\def \de{\delta}

\def \ep{\epsilon}

\def \la{\lambda} 
\def \La{\Lambda}

\def \si{\sigma}

\def \th{\theta}

\def \ff{\infty}
\def \wh{\widehat}
\def \wt{\widetilde}

\def \AA{{\cal A}}
\def \BB{{\cal B}}
\def \CC{{\cal C}}

\def \EE{{\cal E}}
\def \FF{{\cal F}}

\def \NN{{\cal N}}

\def \SS{{\cal S}}

\def \({\left(}
\def \){\right)}

\def \nn{\nonumber}
\def \Proof{\noindent{\bf Proof $\,$ }}

\def \bc{\begin{center} }
\def \ec{\end{center} }
\def \bs{\begin{slide} }
\def \es{\end{slide} }

\def\square{{\vcenter{\vbox{\hrule height.3pt
        \hbox{\vrule width.3pt height5pt \kern5pt
           \vrule width.3pt}
        \hrule height.3pt}}}}
\def\qed{{\hfill $\square$ \bigskip}}

\eqnsection
\begin{document}
  \def\wh{\widehat}
\def\ol{\overline}
\title{ Permanental vectors  }

     \author{\small Hana Kogan\\\small City College of CUNY \and\small  Michael B. Marcus \thanks{Research   supported by  grants from the National Science Foundation and PSCCUNY.}\\\small City College and the  CUNY Graduate Center}
\maketitle

\begin{abstract}    A permanental vector is a generalization of a vector with components that are squares of the components of a Gaussian vector, in the sense that the matrix  that appears in the Laplace transform of the vector of Gaussian squares is not required to be   either symmetric or positive definite. In addition the power of the determinant in the Laplace transform of the vector of Gaussian squares, which is -1/2,  is  allowed to be any number less than zero.

It was not at all clear what vectors are permanental vectors. In this paper we characterize all permanental vectors in $R^{3}_{+}$ and give   applications to permanental vectors in $R^{n}_{+}$ and to the study of permanental processes.  


\end{abstract}

\bibliographystyle{amsplain}

 \section{Introduction}   \label{sec-1}

An  $\bb$-permanental   vector  $\th:=\{\th_{1},\ldots,\th_{n}  \} $,     is an $R_{+}^{n} $  valued random variable with Laplace transform  
 \begin{equation}
   E\(\exp\(- \sum_{i=1}^{n}\al_{i}\th_{ {i}}\)\)=\frac{1}{|I+\al \Gamma |^{\bb}}\label{1.8}, 
      \end{equation}
 where $I$ is the $n\times n$ identity matrix, $\al$ is the diagonal matrix with $  \al_{i,i}=\al_{i} $,  $\al_{i}\in R_{+}$, $1\le i\le n$, and    $\Ga=\{\Ga( {i},  {j})\}_{i,j=1}^{n}$ is an $n\times n$ matrix, $\bb>0$ and $|I+\al \Gamma |>0$ for all   $\al \in R_+^n $.     (This last requirement implies that 
$
   \det \Ga\ge 0$.) The fundamental question here is: For what $n\times n$ matrices is the right-hand side of (\ref{1.8}) a Laplace transform?
   
There are very well known cases in which  right-hand side of (\ref{1.8}) is a Laplace transform. When  $\bb=1/2$ and  $\Ga$ is symmetric and positive definite, 
  \be \th=(G_{1}^{2}/2,  \ldots,G_{n}^{2}/2)
  \ee
 where $(G_{1} ,\ldots,G_{n} )$
   is a Gaussian random variable with covariance $\Ga$. (We sometimes refer to a vector like $\th$ as a vector of Gaussian squares.)

The innovation  in the question posed here is that we consider all $\bb>0$ and do not require that $\Ga$ is symmetric   or positive definite.  

 Before we attempt to answer this question it is important to note that 
the matrix $\Ga$ is not unique.   If $D$ is any diagonal matrix with non-zero entries we have
    \begin{equation}
   |I+\al \Ga |=|I+\al D\Ga D^{-1} |=|I+\al D\Ga^{T} D ^{-1}|,\label{1.3a}
   \end{equation}
   for all diagonal matrices $\al$.  The matrix $ \Ga$ is said to be diagonally equivalent to $\Ga'$ if $\Ga'=D\Ga D^{-1}$ for some diagonal matrix $D$ with non-zero entries.   
  For a very large class of   irreducible matrices $\Ga$, it is known that the class of diagonally equivalent matrices are the only sources of non-uniqueness; see \cite{L}.
  
  Sometimes one can take $D$ to have diagonal entries $\pm 1$. Such matrices are called signature matrices. It is obvious that if $S$ is a signature matrix then $S=S^{-1}$. We also note that (\ref{1.3a}) may hold with $D=I$, the identity matrix even when $\Ga\ne \Ga'$.
   For example, if $\Ga $ and $\Ga'$ are $n\times n$ matrices with the same diagonal elements   and all zeros below the diagonal, then (\ref{1.3a}) holds with $D=I$.     In this case we say that $\Ga $ and $\Ga'$ are effectively equivalent. (We also note that we sometimes refer to $\Ga$ as a kernel for $\th$.)
   
 R. Bapat and R. C. Griffiths \cite{Bapat, Griffiths},  (see, also \cite[Chapter 13]{book}), completely describe the vectors of Gaussian squares for which, (\ref{1.8})   is a Laplace transform for all $\bb>0$. They do this in solving a classical problem posed by P. L\'evy: When is a vector of Gaussian squares infinitely divisible? The answer is:
 
\begin{proposition}\label{prop-BG}	 {\it A vector of Gaussian squares is infinitely divisible if and only if the   covariance matrix $\Ga$ is diagonally equivalent to an $M$-matrix.}
\end{proposition} 

	A matrix
$A=\{a_{i,j}\}$,  is said to be an  
$M$ matrix   if
\begin{enumerate}
\item[(1)] $a_{ i,j}\leq 0$ for all $i\neq j$.
\item[(2)] $A$ is nonsingular and $A^{ -1}\geq 0$.
\end{enumerate}
  
  Strictly speaking knowing that a vector of Gaussian squares is infinitely divisible only asserts the existence of the Laplace transform for $\bb=k/(2n)$, for all integers $k,n\ge 1$.  However the proof of  Proposition \ref{prop-BG} shows that   (\ref{1.8}) holds for all $\bb>0$.
  
 There are permanental vectors with kernels that are not diagonally equivalent to symmetric matrices.  Eisenbaum and Kaspi, \cite[Lemma 4.2]{EK} recognize that the Bapat--Griffiths   sufficient condition for infinite divisibility in the case of symmetric kernels also works for non-symmetric kernels.

\medskip	 	 It is well known that positive definite symmetric matrices may or may not have an inverse that is diagonally equivalent to an $M$-matrix.   On the other hand when $\bb=1/2$, (\ref{1.8}) is the Laplace transform  of a vector of Gaussian squares.  
 Based on these observations we   divide the class of kernels $\Ga$ of permanental vectors into three categories.

\begin{enumerate}
\item $\Ga$ is diagonally equivalent to a  symmetric positive definite matrix.
\item   $\Ga^{-1}$  is diagonally equivalent to an $M$ matrix.  
\item $\Ga$ is not in class 1. or class 2. 
\end{enumerate}

Note that we use the expression positive definite to include what is sometimes called positive semi-definite. Also, we emphasize that classes 1. and 2. are not disjoint.

\medskip	There is ample reason to think that there is an abundance of examples of kernels in class 3. One should be able to take a symmetric positive definite matrix not in class 2. and alter its off diagonal elements very slightly. One then might expect that the altered matrix would be in class 3.  We worked for a long time to find an example of a kernel    of a permanental vector in class 3. but were not successful. We then set out to fully characterize 3$\times 3$ matrices that are kernels    of    permanental vectors. The main result of this paper is the following theorem which states that   for permanental vectors in  $R^{3}_{+}$, class 3. is empty:
 
\begin{theorem} \label{theo-1.1}
A $3\times 3$ matrix that is the kernel    of a permanental vector in $R_{+}^{3}$,   belongs to class 1. or class 2., or both.\end{theorem}

 This result also applies to permanental vectors in $R_{+}^{n}$ in the sense that if  $\th=\{\th_{1},\ldots,\th_{n}  \} $ is a permanental vector  in $R_{+}^{n}$ then any three components of $\th$ is a permanental vector in     $R_{+}^{3}$.

 Another  consequence of 	Theorem \ref{theo-1.1} is that if a kernel of a permanental vector in $R^{3}_{+}$  is not diagonally equivalent to a kernel with positive entries then it is the   kernel   of  a vector of Gaussian squares.
 
\medskip	  A permanental process $\{P(t), t\in T\}$ is a stochastic process with finite dimensional distributions that are permanental vectors.  Eisenbaum and Kaspi study permanental processes in \cite{EK}. Roughly speaking they show that the the potential density of a   Markov process is the kernel of a permanental process. (When this is the case we say that the permanental process is associated with the Markov process.) In fact they show that  permanental processes are the missing link that allows the Dynkin Isomorphism Theorem to be extended to the local times of Markov processes that are not symmetric, \cite[Corollary 3.5]{EK}. There are several intimate connections between permanental processes with a kernel that is the potential density of a Markov process and the Markov process  itself. In \cite{perm} the permanental process is shown to be the loop soup local time of the Markov process.

  
Permanental processes are introduced by 	Vere-Jones in \cite{VJ}.	In \cite[Proposition 4.5]{VJ} he gives necessary and sufficient conditions for (\ref{1.8}) to be the Laplace transform of the vector $(\th _{1} ,\ldots,\th _{n} )$ in terms of the  modified resolvent matrix 
\begin{equation}
   \Ga_{r}:= \Ga(I+r\Ga)^{-1}\label{2.6}
   \end{equation}
where  $r\ge 0$, and $\Ga$ is the matrix in (\ref{1.8}).

\begin{proposition} \cite[Proposition 4.5]{VJ}\label{prop-2.1} For (\ref{1.8}) to represent the Laplace transform of a non-negative random vector it is necessary and sufficient that for all  $r\ge 0$
\begin{itemize}
\item[(i)] $\Ga_{r}$ exists and is $\bb$-positive definite.
\item[(ii)] $ \det(I+r\Ga)>0$.
\end{itemize}
\end{proposition}

Furthermore, given $\Ga_{r}$, Proposition \ref{prop-2.1}, {\em (i)} may hold for some values of $\bb$ but not for others

(Item {\em (ii)} is equivalent to: All the real, non-zero, eigenvalues of $\Ga$ are positive.)

\medskip	There is no point in giving the very complicated definition of $\bb$-positive definite here. One can refer to \cite{VJ} or to \cite{EK,KMJ} where it is repeated. It seems almost impossible to verify  Proposition \ref{prop-2.1} {\em (i)} unless   all the entries of the matrix $\Ga_{r}$ are greater than or equal to zero, in which case  {\em (i)} holds for all $\bb>0$. In \cite[Theorem 3.1]{EK} Eisenbaum and Kaspi point out that this is is the case when $\Ga$ is associated with a Markov process    and that Proposition \ref{prop-2.1}  {\em (ii)}, also holds for these kernels.   

If all the entries of the matrix $\Ga_{r}$ are not greater than or equal to zero   verifying Proposition \ref{prop-2.1} {\em (i)}  necessitates examining  an infinite sequence of increasingly larger matrices derived from  $\Ga_{r}$. Otherwise we know  no sufficient condition for the existence of a permanental vector  that might be in class 3.  One is given in   \cite[Proposition 4.6]{VJ}, but it is not correct. We discuss this in  Remark \ref{VJ_counterex}.

There is a potentially important application of  Proposition \ref{prop-2.1}  if one can figure out how the verify {\em (i)}. If the kernel of a permanental vector is in class 2., (\ref{1.8}) is a Laplace transform for all $\bb>0$. If the kernel of 
a permanental vector is in class 1. and not in class 2. then we only know that (\ref{1.8}) is a Laplace transform for   $\bb=1/2$, and  trivially, for all $\bb=k/2$, for integers $k\ge 1$.  Possibly there exist  other values of $\bb>0$ for which (\ref{1.8}) is a Laplace transform. Applying Proposition \ref{prop-2.1}, which depends on $\bb$, would answer this question.

\medskip		
	There are many other interesting applications of  Theorem \ref{theo-1.1}. The next result answers a question that  started our interest in 3-dimensional permanental vectors. We point out in the beginning of this Introduction  that the univariate marginals 
of a $1/2$-permanental process are squares of normal random variables.  It also follows  from (\ref{1.8}) that  pairs $(\th_{i},\th_{j})$, of a $1/2$-permanental process, are equal in law to $(G_{i}^{2}/2,G_{j}^{2}/2)$, where  $(G_{i},G_{j})$ is a Gaussian vector with covariance matrix
 \be\wt\Ga =\begin{bmatrix}
\Ga(i,i)& \(\Ga(i,j)\Ga(j,i)\)^{1/2}\\ 
\(\Ga(i,j)\Ga(j,i)\)^{1/2}&\Ga(j,j) 
\end{bmatrix}\label{1.3}\ee
 (See \cite[Lemma 3.1]{perm}). It follows from this that 
   \begin{equation}
 E(\th _{ i})=\frac{\Ga (i,i)}{2} \qquad  \mbox{and}\qquad\mbox{cov}\{\th _{ i},\th_{j}\}=\frac{ \Ga(i,j)\Ga(j,i)}2.\label{4.2v}
   \end{equation}
   Therefore, if 
   \begin{equation}
   \Ga(i,j)\Ga(j,i)=0\qquad \forall \,1\le i\ne j\le n\label{1.4}
   \end{equation}
   the components of  a $1/2$-permanental process are pairwise independent. 
  
  Actually,  (\ref{1.4})  is a necessary and sufficient condition for  the components of  any $\bb$-permanental process to be pairwise independent. This is because in this case the determinant of  $\wt \Ga$
is   a product of its diagonal elements and the right-hand side of (\ref{1.8})
\begin{equation}
   |I+\al\Ga|^{\bb}= |I+\al_{i}\Ga(i,i)|^{\bb} |I+\al_{j}\Ga(j,j)|^{\bb}.
   \end{equation}
  We also know from  \cite[bottom of page 135]{VJ} that for any $\bb$ permanental process
   \begin{equation}
 \qquad\mbox{cov}\{\th _{ i},\th_{j}\}=\bb\Ga(i,j)\Ga(j,i).\label{4.2vn}
   \end{equation}

      If $\Ga$ is symmetric and positive definite and $\th=(G_{1}^{2}/2,\ldots,G_{n}^{2}/2)$, where $(G_{1} ,\ldots,G_{n} )$ is a Gaussian random variable with covariance $\Ga$, with $   \Ga(i,j)=\Ga(j,i)=0$, then   the components of  $\th$ are independent. We asked ourselves the following question: ``For a general $\bb$-permanental vector $\th$, that is not the square of a Gaussian vector,  does  (\ref{1.4}) imply that the components of  $\th$ are independent?'' The answer is yes. We prove:
   
   \begin{theorem}\label{theo-1.2}
Let $\th$ be an $n$-dimensional $\bb$-permanental vector with pairwise independent components. Then the components of $\th$ are independent.
\end {theorem}

 It is clear that  when $\bb=1/2$ Theorem \ref{theo-1.2} implies that the only permanental vectors with independent components are those with components that are   squares of independent Gaussian random variables.

\medskip	
 The next result deals with a function that appears in sufficient conditions for the continuity of permanental processes in \cite{perm},
  \begin{equation}
   d(x,y) = \(\Ga(x,x)+\Ga(y,y)-2\(\Ga(x , y)\Ga(y, x )\)^{1/2}\)^{1/2}.\label{1.2w}
   \end{equation}
If $\Ga(x , y)=\Ga(y, x )$ is the covariance of the Gaussian vector $\{G(t), t\in T\}$, then
\begin{equation}
    d(x,y)=\(E\(G(x)-G(y)\) ^{2}\)^{1/2},\label{1.7}
   \end{equation}
which is a metric on $T$. However, if $\Ga(x , y)$ is the kernel of a permanental  process and $ \Ga(x , y)\ne \Ga(y, x )$ it was not clear whether 
or not $   d(x,y)$ is a metric on $T$. We can now say that even in this case $   d(x,y)$ is a metric on $T$.

\begin{corollary}\label{cor-1.1} Let  $ \{P(t), t\in T\}$ be a permanental   process  with kernel $\Ga(x , y)$. The function $\{ d(x,y), x,y\in T\}$ in  (\ref{1.2w}) is a metric on $T$.

 \end{corollary}

 In Section  \ref{sec-2} we give many properties that are necessary for an $n\times n$ matrix to be the kernel of a permanental vector.   In Section  \ref{sec-3} we obtain an interesting property of the eigenvalues of  ${ 3\times 3}$ positive definite symmetric matrices that plays a critical role in the proof of Theorem \ref{theo-1.1}. The proof of Theorem \ref{theo-1.1} uses completely different methods when the off diagonal elements of the kernel are all negative or all positive. These cases are considered separately in  Sections  \ref{sec-4} and \ref{sec-5}. Sections \ref{sec-6} and \ref{sec-7} give, repectively,  the proofs of  Theorem \ref{theo-1.2} and  Corollary  \ref{cor-1.1}. Finally, because the fact that kernels in class 2. are kernels of permanental vectors is so important in this paper, we give an  outline the proof, essentially showing what changes are necessary    in the proof in \cite[Theorem ]{book}, which is given for symmetric kernels.   This result is given in \cite[Lemma 4.2]{EK}.   The  proof involves probabilistic considerations. Since Theorem \ref{theo-1.1} is only for $3\times 3$ matrices it seems appropriate to give a proof for finite matrices involving only linear algebra.  
 
 \medskip	We are grateful to Professor Jay Rosen for many helpful comments and discussions

\section{Preliminaries}\label{sec-2}

If $\th$ is a permanental vector in $R^{n}_{+}$   then   any subset  of its components, say of $p$ components, is a permanental vector in $R^{p}_{+}$.
For $p=2$, the Laplace transform of the vector $\{\th_{i},\th_{j}\}$  takes the form  
  \bea
  && E\(\exp\(-\frac{1}{2}\(\al_{i} \th_{i }+\al_{j} \th_{j}\)\)\)\label{i.8a}\\
  &&\qquad=\frac{ 1}{|I+\al \Gamma |^{\bb}}\nn= \(1+\al_{i}\Ga(i,i))+\al_{j}\Ga(j,j)\right.\\
& &   \nn \qquad\qquad\left.+\al_{i}\al_{j}\(\Ga(i,i )\Ga(j,j)-\Ga(i,j)\Ga(j,i )\)\)^{-\bb}.
   \eea
   
   Taking $\al_{i}=\al_{j}$  sufficiently large, this implies that
    \be 
\Ga(i,i)\Ga(j,j)-\Ga(i,j)\Ga(j,i)\geq 0. \label{3.2k}
 \ee
If we also set  $\al_{j}=0$ in (\ref{1.8}) we see that for any   $i\in n$
 \begin{equation}
 \Ga(i,i )\geq 0\label{vj.1}.
 \end{equation}
 In addition, by \cite[Proposition 3.8]{VJ},  for any pair $i,j\in T$  
  \begin{equation}
 \Ga(i,j ) \Ga(j,i )\geq 0.\label{vj.2}
 \end{equation}

In the next lemma  we show that there are many transformations of kernels of permanental processes that give other  kernels of permanental processes. 
 
\bl \label{lem-2.3}Let $A$ be a kernel of a $\bb$-permanental vector $\th=(\th_{1},\ldots ,\newline \th_{n})$. Let $U_{1}$ and $U_{2}$ be     diagonal matrices with non-zero diagonal entries $u^{(j) }_{i }$, $i=1,\ldots, n$, $j=1,2$, for $U_{1}$ and $U_{2}$ respectively,  with the property that $u^{(1)}_{i }u^{(2)}_{i }>0$, $i=1,\ldots, n$. Then     $U_{1}AU_{2}$ is the kernel of the $\bb$-permanental vector $ (u^{(1)}_{1 }u^{(2)}_{1 } \th_{1},\ldots ,u^{(1)}_{n }u^{(2)}_{n}\th_{n})$. 
 \el

 \Proof  Since $\th$ is an $R^{n}_{+}$ valued random variable so is $ (u^{(1)}_{1 }u^{(2)}_{1 } \th_{1},\ldots ,\newline u^{(1)}_{n }u^{(2)}_{n}\th_{n})$.
The Laplace transform of $ (u^{(1)}_{1 }u^{(2)}_{1 } \th_{1},\ldots ,u^{(1)}_{n }u^{(2)}_{n}\th_{n})$ is
\bea
E\(\exp\(-\sum_{i=1}^{n} \al_{i} ( u^{(1)}_{i }u^{(2)}_{i }\th_{i})\)\)&=&E\(\exp\(-\sum_{i=1}^{n}( \al_{i} u^{(1)}_{i }u^{(2)}_{i })\th_{i})\)\)\nn\\
&=& |I+(\al U_{2}U_{1} )A|^{ -\beta} \\
&=&|U_{ 2}(I+\al (U_{ 1}A U_{ 2})U_{2}^{-1}|^{ -\beta}\nn\\
&=&| I+\al (U_{ 1}A U_{ 2}) |^{ -\beta}\nn.
\eea
\qed

\begin{example}\label{ex-2.1} {\rm We note two cases. Let   $U  $ be   a strictly positive diagonal matrix
\begin{enumerate}
\item When $\ga+\ga'$=0 \begin{equation} E\(\exp\(-\sum_{i=1}^{n} \al_{i}  \th_{i}\)\)=| I+\al  U^{\ga}A U^{-\ga} |^{ -\beta}.
   \end{equation}
   \item When $\ga+\ga'$=1 \begin{equation} E\(\exp\(-\sum_{i=1}^{n} \al_{i} u_{i} \th_{i}\)\)=| I+\al  U^{\ga}A U^{ (1-\ga)} |^{ -\beta}.
   \end{equation}
   In particular,
    \begin{equation} E\(\exp\(-\sum_{i=1}^{n} \al_{i} u_{i} \th_{i}\)\)=| I+\al  U^{1/2}A U^{  1/2} |^{ -\beta}.
   \end{equation}

 \end{enumerate}
 }
 \end{example}
 
\begin{remark}\label{ex-2.2} {\rm  It is easy to see that for $a',b',c'$ strictly positive, the two   matrices  
 \begin{equation}
\(\begin {array}{ccc}
1&a'&c'\\
a'&1&b'\\
c'&b'&1
\end{array}\) \label{3n.11}  \mbox{and} \(\begin {array}{ccc}
1&-a'&-c'\\
-a'&1&b'\\
-c'&b'&1
\end{array}\).
\ee
  are  diagonally equivalent to each other. Similarly  
\be
\(\begin {array}{ccc}
1&-a'& c'\\
-a'&1&b'\\
 c'&b'&1
\end{array}\)  \mbox{and}   \(\begin {array}{ccc}
1&-a'&-c'\\
-a'&1&-b'\\
-c'&-b'&1
\end{array}\)\label{3n.12}
\end{equation}
  are  diagonally equivalent to each other. It should also be clear that these observations hold if any two of the three pairs of  entries in (\ref{3n.11}) are taken to be negative and if any   pair  of  entries in (\ref{3n.12}) is taken to be negative. 
  
  Because of the observations in the previous paragraph, when we consider whether a $3\times 3$ matrix is the kernel of a permanental vector we need only consider those matrices with all positive off-diagonal elements or all negative off-diagonal elements. (We consider 0 to be both positive and negative.)
 }\end{remark}
 
 Consider the matrix

\begin{equation}
E= 
\(\begin {array}{ccc}
1&a_1&c_2\\
a_2&1&b_1\\
c_1&b_2&1
\end{array}\).\label{2.7}
\end{equation}

 The next elementary lemma is very useful. We leave the proof to the reader.   \begin{lemma} \label{lem-2.3a}If the off diagonal elements of the  matrix $E$ in (\ref{2.7}) are either all strictly positive or all strictly negative   then $E$ is diagonally equivalent to
\begin{equation}
E'= 
\(\begin {array}{ccc}
1&a &c \\
a &1&b'_1\\
c &b'_2&1
\end{array}\).\label{2.7hh}
\end{equation}
where $a^{2}=a_{1}a_{2}$, $c^{2}=c_{1}c_{2}$ and $b'_{1}b'_{2}=b _{1}b _{2}$ and the signs of $a,c,b'_1,b'_2$ are the same as the signs of $a_{1},a_{2}, c_{1},c_{2}, b _{1},b _{2}$.
 \end{lemma}

 By (\ref{vj.2})   the kernel $\Ga$ of a permanental vector   has the property that $\Ga(i,j)\Ga(j,i)\ge 0$. Therefore, if $E$ is the kernel  of a permanental vector in $R_{+}^{3}$,      $a_{1},a_{2}$ are either both positive or both negative, and similarly for $b_{1},b_{2}$ and $c_{1},c_{2}$.

\begin{lemma} \label{lem-7.1} Suppose that the matrix $E$ is the kernel of a permanental vector in $R_{+}^{3}$. Then if
\begin{equation}
   a_{1}b_{1}c_{1}=a_{2}b_{2}c_{2}\label{7.3}
   \end{equation}
   it is  diagonally  equivalent to the kernel   
 \begin{equation}
 \EE= 
\(\begin {array}{ccc} 
1&\pm a&\pm c\\
\pm a&1&\pm b\\
\pm c&\pm b&1
\end{array}\) 
\end{equation} 
where $a=(a_1a_2)^{1/2}$, $b=(b_1b_2)^{1/2}$ and $ c=(c_1c_2)^{1/2}$ and  
in which  $\{\EE\}_{1,2}=a$  if $a_{1}$ is positive and $\{\EE\}_{1,2}=-a$ if $a_{1}$ is negative, and similarly with respect to  $b_{1} $ and $c_{1} $.

 In particular this lemma  holds when both sides of (\ref{7.3}) are equal to zero. 
 \end{lemma}
 
 Obviously,  $\EE$ is the covariance of a Gaussian vector. 
 
 \medskip	 
 
 \Proof When (\ref{7.3}) holds and  $a_{1}b_{1}c_{1}\ne 0$ 
 \begin{equation}
a_{1}b_{1}c_{1}+a_{2}b_{2}c_{2}=2(a_{1}a_{2}b_{1}b_{2}c_{1}c_{2})^{1/2}\label{2.14}
   \end{equation}
    It is easy to see that  $E$ and $\EE$ are diagonally equivalent.
  It is also easy to see that   $E$ and $\EE$ are diagonally equivalent 
 if, say, $b_{1}=b_{2}=0$.
 
  Finally,  it is also easy to see that if $a_{1}=b_{2}=0$, $ E$ is effectively equivalent to \be \EE'= 
\(\begin {array}{ccc}
1&  0&  (c_1c_{2}) ^{1/2}\\
0&1&0 \\
  (c_1c_{2}) ^{1/2}& 0&1
\end{array}\)\label{3.1aw}.
\end{equation} 
 \qed

 We also use the following lemma which is   \cite[Lemma 4.5]{KMJ}.

 \bl \label{lem-1.1}Let 
\be  \AA=\left (   \begin{array}{cccc}
    u& a& c\\ 
   a& v& b  \\ 
 c&b & w \\ 
  \end{array} \right ),\qquad \BB=\left (   \begin{array}{cccc}
    u&a_{1}&c_{2}\\ 
 a_{2}& v& b_{1} \\ 
c_{1}&b_{2} & w \\ 
  \end{array} \right ).
\ee
where   $a_{1}a_{2}=a^{2}$,   $b_{1}b_{2}=b^{2}$,$c_{1}c_{2}=c^{2}$. Suppose that $\AA\ge 0$. If $\BB^{-1}$ is an $M$ matrix then $\AA^{-1}$ is an $M$ matrix.
\el

\begin{remark}\label{rem-2.2} {\rm 
By definition an $M$-matrix is invertible. Therefore   $|\AA|>0$. We also note that if $\BB^{-1}$ is an $M$ matrix then for any diagonal matrix $D$ with strictly positive entries, $(D\BB D^{-1})^{-1}$  is an $M$-matrix.
 }\end{remark}

 \medskip	
 The next  observation is used often in this paper.

 \begin{lemma}\label{lem-2.1} Let $\Phi(\al_{1},\ldots,\al_{n})$ be the Laplace transform of an $R_{+}^{n}$ valued random variable. For any $1<k< n$ set $\al_{j}=u_{j}$,  where $ u_{j}\ge 0$, $k\le j\le n$. Then 
  \be
   \Phi_{(n,k)}(\al_{1},\ldots,\al_{k }) =\frac{\Phi(\al_{1},\ldots,\al_{k },u_{k+1},\ldots,u_{n})}{ 
 \Phi(0,\ldots,0,u_{k+1},\ldots,u_{n}) }\label{2.5}
 \ee
is the Laplace transform of an $R_{+}^{k}$ valued random variable. 

  Furthermore, if $\Phi(\al_{1},\ldots,\al_{n})$ is the Laplace transform of an $n$-dimen-\newline sional permanental vector, $ \Phi_{(n,k)}(\al_{1},\ldots,\al_{k }) $  is the Laplace transform of a $k$-dimensional permanental vector.
 \end{lemma}
 
 \Proof Since $\Phi(\al_{1},\ldots,\al_{n})$ is a completely monotone function on $R_{+}^{n}$ it follows that $ \Phi_{(n,k)}(\al_{1},\ldots,\al_{k })$
  is a completely monotone function on $R_{+}^{k }$, satisfying $ \Phi_{(n,k)}(0,\ldots,0)=1$. Therefore, it is the Laplace transform of an $R_{+}^{k}$ valued random variable.
  
     (This is very well known when $k=1$. Lacking a suitable reference for general $k$, we note that it follows from the Extended Continuity Theorem for probability measures on $R^k_+$,  \cite[Theorem  5.22]{Ka}, and the argument in the proof  of \cite[Lemma  13.2.2]{book}, applied to $ \Phi_{(n,k)}(\al_{1},\ldots,\al_{k })$, not its logarithm.)
 
 Now suppose that $\Phi(\al_{1},\ldots,\al_{n})$ is the Laplace transform of an $n$-dimensional permanental vector. This implies that
 \begin{equation}
   \Phi(\al_{1},\ldots,\al_{n})=\frac{1}{|I+\al \Ga|^{\bb}}\label{2.19}
   \end{equation} 
  for some $n\times n$ matrix $\Ga$, and diagonal matrix $\al$ as in (\ref{1.8}). We first prove the second statement in the lemma for $k=n-1$. Consider 
  \begin{equation}
   \Phi  (\al_{1},\ldots,\al_{n-1 },u_{n}) 
     \end{equation} 
   and the corresponding matrix $I+\wt \al \Ga$, where $\wt \al=(\al_{1},\ldots,\al_{n-1},u_{n})$. 
   
 	We now show that 
   \begin{equation}
  | I+\wt \al \Ga |=(1+{u_{n}}) | I+  \al^{(n-1) }\Ga^{(n-1) } |\label{2.8}
   \end{equation}
   where $ \al^{(n-1)}$ is the $(n-1)\times (n-1)$ diagonal matrix with diagonal entries $(a_{1},\ldots,\al_{n-1})$ and $\Ga^{(n-1) }$ is an $(n-1)\times (n-1)$   matrix with entries that are functions of the entries of $\Ga$ and $u_{n}$. Since
  \begin{equation}
    \Phi(0,\ldots,0, \ldots,0,u_{n}) =\frac{1}{|1+u_{n}|^{\bb}}
   \end{equation}
  The equality in  (\ref{2.8}) gives (\ref{2.5}) when $k=n-1$.
   
  To obtain (\ref{2.8}) we note the matrix $I+\wt \al \Ga$ has the same determinant as the matrix obtained from it by subtracting   $\Ga(n,j) \frac{u_{n}}{(1+u_{n})}$ times the $n$-th column from the $j$-th column,   for each $1\leq j \leq n-1$. Call this matrix $S$. Note that $S(n,j)=0$, $j=1,\ldots,n-1$ and $S(n,n)=1+u_{n}$. Let $S'$ denote the martix obtained by dividing the last row of $S$ by $ 1+u_{n}$.  We have
  \begin{equation}
    | I+\wt \al \Ga |=(1+u_{n})   | S' |.
   \end{equation}
  
   To be more specific the entries of $S'$ are
   \bea
    S'(i,j)&=& \de_{i,j}+\al_i \(\Ga( i, j)-\frac{u_{n}\Ga(i,n )\Ga(n,j)}{1+u_{n}}\),\quad 1\le i,j\le n-1;\nn\\
S'(n,j)&=&0,\qquad\quad\hspace{2.2in} 1\le  j\le n-1;\nn\\
       S'(n,n)&=&1\label{2.10}.
         \eea
 It is obvious that we can write   
 \bea 
  |S'| =| I+  \al^{(n-1) }\Ga^{(n-1) } |
 \eea 
  where $\Ga^{(n-1) }$ is the matrix with components
   \be 
  \(\Ga( i, j)-\frac{u_{n}\Ga(i,n )\Ga(n,j)}{1+u_{n}}\),\qquad 1\le i,j\le n-1\label{2.26-}.
   \ee 
 and $\al^{(n-1)}=( \al_1, \al_2, ...\al_{n-1}).$ 
   We now have 
   \begin{equation}
     \Phi_{(n,n-1)}(\al_{1},\ldots,\al_{n-1 }) =\frac{1}{| I+  \al^{(n-1) }\Ga^{(n-1) } |^{\bb}}.\label{2.26}
   \end{equation}
   Repeating the argument above we can show that 
   \bea 
      \Phi_{(n,n-2)}(\al_{1},\ldots,\al_{n-2 })&=&\frac{   \Phi_{(n,n-1)}(\al_{1},\ldots,\al_{n-2 },u_{n-1})}{   \Phi_{(n,n-1)}(0,\ldots,0,u_{n-1 })}\\
      &=&\frac{   \Phi (\al_{1},\ldots,\al_{n-2 },u_{n-1},u_{n})}{   \Phi (0,\ldots,0,u_{n-1 },u_{n})}\nn,
 \eea
  since 
  \begin{equation}
  \Phi_{(n,n-1)}(\al_{1},\ldots,\al_{n-2},u_{n-1}) =\frac{\Phi(\al_{1},\ldots,\al_{n-2}, u_{n-1},u_{n})}{ 
 \Phi(0,\ldots, 0 ,u_{n}) }\label{2.5qq}
   \end{equation}
   and
   \begin{equation}
  \Phi_{(n,n-1)}(0,\ldots,0,u_{n-1 })=\frac{\Phi(0,\ldots,0, u_{n-1},u_{n})}{ 
 \Phi(0,\ldots, 0 ,u_{n}) } .
   \end{equation}
  Thus we get (\ref{2.5}) for $k=n-2$. 
  
  Continuing in this way we get  (\ref{2.5}) for arbitrary $1\le k\le n-1$. \qed

  We use the following necessary condition in the proof of Theorem \ref{theo-1.1}.  It is a direct consequence of \cite[Proposition 3.8 and Proposition  4.5  with $\si$=0]{VJ}.  We provide a direct proof for the convenience of the reader.

  \bl  \label{lem-2.2}
Let  $A=\{A_{i,j} \}_{i,j=1}^{n}$, be an $n \times n$ matrix. If $A$  is a kernel of a $\bb$-permanental vector then $A$ and all matrices obtained from $A$ by multiplying its rows by non-negative numbers have a  positive  eigenvalue of maximum modulus. \el

\Proof Let $\th=(\th_{1}, \dots,\th_{n})$ be a permanental vector with kernel $A$. The matrices obtained by multiplying the rows of $A$  by non-negative numbers have the form $UA$ where $U$ is a diagonal matrix with non-negative entries $u_{1}\ldots,u_{n}$.  Note that
\begin{equation} E\(\exp\(-\sum_{i=1}^{n} \al_{i} u_{i} \th_{i}\)\)=| I+\al  U A   |^{ -\beta}.\label{4.1a}
   \end{equation}
  Let $z$ be a complex number and set 
\bea
f(z)&=&E\(\exp\( {z} \sum_{i=1}^{n}   u_{i} \th_{i}\)\) \label{4.1}\\&=& \sum_{k=0}^\infty  z^k\frac{E\(\sum_{i=1}^{n}   u_{i} \th_{i}\)^k}{k!}. \nn
\eea
By (\ref{4.1a}) 
\be \label{2}
f(z)  = |I- zU A|^{ \beta}= \prod_{p=1}^n {(1-z \lambda_p)^{- \beta}} 
\ee
where   $\la _p$, $1\le p\le n$,  are the eigenvalues of $U A$.

Since  
    \begin{equation}
   E\(\exp\(- \la\th_{ {i}} \)\)=\frac{1}{|I+\la A_{i,i}|^{\bb}}\label{ }, 
      \end{equation}
 we  see that
 \begin{equation}
   E(\th_{x_{i}}^{k})=(A_{i,i}^k)\Ga( \beta+k).
   \end{equation}
It follows that
\bea
 E\(\sum_{i=1}^n {u_i \th_{i}}\)^k&\leq& \max_{1\le i\le n} u_i ^k n^k \max_{1\le i\le n}E(\th_{ {i}}^k) \\
& =& \max_{1\le i\le n}u_i^k n^k \max_{1\le i\le n}A_{i,i}^k\Ga( \beta+k).\nn
\eea
Consequently, there exists a number b, such that  
  \begin{equation}
 \BB_{k}:={  E\(\sum_{i=1}^n {u_i \th_{i}}\)^k\over k!}\le b^{k}.
   \end{equation}
  This implies that the series in (\ref{4.1}) has a positive radius of convergence, which we denote by $R$. By (\ref{4.1}) and (\ref{2})
\be
f(z)=\sum _{k=0}^\infty  \BB_{k} z^k =\prod_{p=1}^n {(1-z \lambda_p)^{- \beta}} \label{3}
\ee
for $|z|<R$. Let $v>0$ and  note that $\lim_{v \rightarrow R} f(v) = \ff$ as the sum of a series with positive terms.

 Since the terms $\BB_k$ are positive, when $v=|z|$ we have  
   \be\label{4}
 f(v)=   \sum _{k=0}^\ff  \BB_{k}  |z| ^k \geq   \left | \sum _{k=0}^\ff  \BB_{k}  z ^k \right|.
   \ee 
This shows that if    $f(z)$ has a singularity at $z_0$, then $f(v_{0})=\ff$, for  $v_{ 0}=|z_{0}|$.

 By (\ref{3}) this implies that 
\begin{equation}
  \max_{1\le p\le n}\la_{p}=\frac{1}{v_{0}}.
   \end{equation}\qed

\section{Eigenvalues of  ${\bf 3\times 3}$ positive definite symmetric matrices}\label{sec-3}

\bl \label{general}
Let $     H$   be the  real symmetric matrix  
\begin{equation}
  H= 
\(\begin {array}{ccc}
1&a&c\\
a&1&b\\
c&b&1
\end{array}\) \label{3.1}
\end{equation}
  and let $\rho$ be a real diagonal matrix with entries $(\rho_1, \rho_2, \rho_3)$, where
 \be
 \rho_1=\frac{b}{b-ac}, \quad \rho_2=\frac{c}{c-ab},  \quad\mbox{and}\quad\rho_3=\frac{a}{a-bc}.\label{7.19aa}
 \ee
  Assume none of the   denominators in (\ref{7.19aa}) are zero.  Then
 \begin{equation}
 |\rho H-\la I|=(\la - 1)^2(\rho_1\rho_2\rho_3 |H| -\la ).\label{7.20}
      \end{equation}
In particular, $\la=1$  is an eigenvalue of  $\rho H$ of multiplicity 2. 
 \el
  
 \Proof  
 \bea 
 |\rho H-I|= \left|\begin {array}{ccc}
\frac{ac}{b-ac}&\frac{ab}{b-ac}&\frac{bc}{b-ac}\\
\frac{ac}{c-ab}&\frac{ab}{c-ab}&\frac{bc}{c-ab}\\
\frac{ac}{a-bc}&\frac{ab}{a-bc}&\frac{bc}{a-bc}
\end{array}\right| .\eea
It is easy to see that the second and third row of this determinant are equal to a (different) multiple of the first row. This shows that $\la=1$  is an eigenvalue of  $\rho H$ of multiplicity 2. 

Since the product of the eigenvalues must equal $|\rho H|$ we get (\ref{7.20}).\qed

  \begin{remark} {\rm 
When $\det H\ge 0$, or  equivalently,  when $H$ is positive definite, $H$ is the covariance of a Gaussian vector, say, $(\xi_{1},\xi_{2},\xi_{2})$.  When  $\rho\ge 0$, it follows from Lemma \ref{lem-2.3} that $\rho H$ is diagonally equivalent to the covariance matrix of the Gaussian vector \begin{equation}
   \(\rho_{1} ^{1/2}G_{1},\rho_{2} ^{1/2}G_{2},G_{3} ^{1/2}\xi_{3}\).
   \end{equation}

 }\end{remark}

\section{$\bf3\times 3$ matrices with negative off diagonal elements}\label{sec-4}

In this section we prove Theorem \ref{theo-1.1} for matrices with negative off diagonal elements.  

In the next lemma we consider the eigenvalues of a $3\times 3$ matrix with negative off diagonal elements.  
 
 \bl \label{general1}
Let $ A_-$ be the   matrix  
\begin{equation}
 A_-= 
\(\begin {array}{ccc}
1&-a'_1&-c'_2\\
-a'_2&1&-b'_1\\
-c'_1&-b'_2&1
\end{array}\).\label{2.7w}
\end{equation}
with  $a_{1}'a_{2}',b_{1}'b_{2}',c_{1}'c_{2}'$ all greater than or equal to 0 and less than or equal to 1. Assume that $\det A\ge 0$ and  $A_{i,j}A_{j,i}\leq 1$. 
Then if $ A_-$ is not diagonally equivalent to a symmetric matrix there exists a  diagonal matrix $\Phi$, with strictly positive  entries, such that $\Phi A_-$ has only one real eigenvalue. 

   Furthermore,  the real part of the complex eigenvalues of $\,\Phi A_-$ is greater than the real eigenvalue. 
\el

\Proof  We first consider the case in which $a_{1}',a_{2}', b_{1}',b_{2}',c_{1}',c_{2}'$ are all strictly positive, Let  $a_{1}'a_{2}'=a^{2}$, $b_{1}'b_{2}'=b^{2}$ and $c_{1}'c_{2}'=c^{2}$ and let $ \Phi$ be a real diagonal matrix with entries $(  \phi_1,   \phi_2,   \phi_3)$, where
 \be
  \phi_1=\frac{b}{b+ac}, \quad   \phi_2=\frac{c}{c+ab}  \quad\mbox{and}\quad  \phi_3=\frac{a}{a+bc}.\label{7.19aab}
 \ee
 Let $K$ be the matrix
 \begin{equation}
 K= 
\(\begin {array}{ccc}
1&-a&-c\\
-a&1&-b\\
-c&-b&1
\end{array}\). \label{3.1a}
\end{equation}
 We show in Lemma \ref{general} that
 \bea
\lefteqn{|\Phi K-\la I|}\label{7.20ss}\\
&=&   -\la^3+\la^2\( \phi_1+ \phi_2+ \phi_3 \)\nn  \\  
& &\hspace{.2in} -\la \(  \phi_1 \phi_2(1-a^2)+ \phi_2 \phi_3(1-b^2)+ \phi_1 \phi_3(1-c^2)\)+ \phi_1 \phi_2 \phi_3 |K|\nn\\
 &=&(\la - 1)^2(\phi_1\phi_2\phi_3  |K|-\la ).\nn
      \eea
   Since $(\phi_1+ \phi_2+ \phi_3)/3<1$, the second derivative of  $|\Phi K-\la I|$ is negative at $\la=1$. Therefore, $|\Phi K-\la I|$ has a local maximum of 0 at $\la=1$.
 Consider 
 \be
H(\la):= |\Phi  A_--\la I|.
 \ee
By Lemma  \ref{lem-2.3a} it suffices to take $ A_-$ to be 
\begin{equation}
\(\begin {array}{ccc}
1&-a &-c \\
-a &1&-b_1\\
-c &-b _2&1
\end{array}\).\label{2.7hh1}
\end{equation}
where   $b _{1}b _{2}=b ^{ 2}$. 

  We have
\be 
 |K|=1-a^{2}-b^{2}-c^{2}-2abc 
 \ee 
 and
  \be 
|  A_-|=1-a^{2}-b^{2}-c^{2}- ac (b_{1}+b_{2}).
 \ee 
Therefore, since  $b _{1}b _{2}=b ^{ 2}$ we see from the first equality in (\ref{7.20ss}) that    
\begin{equation}
  |\Phi K-\la I|-|\Phi  A_--\la I|=\phi_1\phi_2\phi_3\(ac(b_{1}+b_{2}-2b)\).\label{3.30}
   \end{equation}
 Unless, $b_{1}=b_{2}$, the right-hand side of (\ref{3.30}) is strictly positive. Since $|\Phi K-\la I|$ has a local maximum of 0 at $\la=1$,  $|\Phi  A_--\la I|$ has a local maximum that is strictly negative at $\la=1$. This implies that when $b_{1}\ne b_{2}$, $|\Phi  A_--\la I|$ has only one real root. This is equivalent to the statement of this lemma.  
 
  We now consider the cases in which some of the entries of $ A_-$ in (\ref{2.7w})
 are zero. Considering Lemma  \ref{lem-2.3a} we can restrict our attention to the following matrices
 \begin{equation}
\(\begin {array}{ccc}
1&-a &-c \\
-a &1&-b' \\
-c &0&1
\end{array}\) \label{2.7dd}\quad   
\(\begin {array}{ccc}
1&-a &0 \\
-a &1&-b' \\
-c' &0&1
\end{array}\) \quad\(\begin {array}{ccc}
1&-a' &0 \\
0 &1&-b' \\
-c' &0&1
\end{array}\). 
 \ee
in which $a,b,c,a',b' ,c'$ are all strictly positive . Label them, respectively, $\FF_{1}$, $\FF_{2}$ and  $\FF_{3}$. 

We first show that $\FF_{1}$ has only one real eigenvalue.  Similar to  (\ref{3.30}) we have
 \bea
&& |\Phi  K-\la I|- |\Phi   \FF_{1}-\la I|\label{7.20aa}\\
&&\qquad=\nn -\phi_{2}\phi_{3}b^{2}\(\phi_{1}-\la\)-\phi_{1}\phi_{2}\phi_{3} \( 2abc-acb'\)\\
&&\qquad=\nn  \phi_{2}\phi_{3}b^{2} \la -\phi_{1}\phi_{2}\phi_{3} \( b^{2}+2abc-acb'\).
      \eea
  We assume $b'$ is fixed.  One can choose $0<b<1$ such that $b^{2}+2abc-acb'=0$. That $b>0$ is elementary. That $b<1$ follows from the fact that $\det \FF_{1}\ge 0$ implies that $acb'<1$. With this choice of $b$
 \begin{equation}
|\Phi  K-\la I|- |\Phi   \FF_{1}-\la I|= \phi_{2}\phi_{3}b^{2} \la.\label{sds}
   \end{equation}
    Considering the graph of $|\Phi K-\la I|$; (see the last line of (\ref{7.20ss})),  and the fact that the right-hand side of   (\ref{sds}) is strictly positive for $\la>0$      we see that $|\Phi   \FF_{1}-\la I|$ has only one real root, or equivalently, that $\FF_{1}$ has only one real eigenvalue. 
    
    Similar arguments show that $\FF_{2}$ has only one real eigenvalue. With regard to $\FF_{2}$ we have
    \bea
&& |\Phi  K-\la I|- |\Phi  \FF_{2}-\la I|\label{7.20aa}\\
&&\qquad=\nn -\phi_{2}\phi_{3}b^{2}\(\phi_{1}-\la\)-\phi_{1}\phi_{3}c^{2}\(\phi_{2}-\la\)-\phi_{1}\phi_{2}\phi_{3} \( 2abc-ac'b'\)\\
&&\qquad=\nn  (\phi_{2}\phi_{3}b^{2}+\phi_{1}\phi_{3}c^{2}) \la -\phi_{1}\phi_{2}\phi_{3} \( b^{2}+c^{2}+2abc-ac'b'\).
      \eea
  For simplicity we can take $b=c$. We can can choose $1>b>0$ such that $2b^{2}+2ab^{2}-ac'b'=0$. With this choice of $b$
   \begin{equation}
|\Phi  K-\la I|- |\Phi  \FF_{2}-\la I|=(\phi_{2}\phi_{3}+\phi_{1}\phi_{3})b^{2} \la.\label{sds1}
   \end{equation}
   So we can use the same argument we just used to show that    $\FF_{2}$ has only one real eigenvalue. 

For  $\FF_{3}$ we consider
  \bea
&& |\Phi  \FF_{2}-\la I|- |\Phi  \FF_{3}-\la I|\label{7.20mm}\\
&&\qquad=\nn  \phi_{1}\phi_{2}a^{2}  \la -\phi_{1}\phi_{2}\phi_{3} \( a^{2} +ab'c'-a'b'c'\).
      \eea
  We can choose $0<a<1$   so that $a^{2} +2ab'c'-a'b'c'=0$ and get 
   \bea
  |\Phi   \FF_{2}-\la I|- |\Phi   \FF_{3}-\la I|  =   \phi_{1}\phi_{2}a^{2}  \la .       \eea
     Considering the graph of $  |\Phi   \FF_{2}-\la I|$ we see that  $\FF_{3}$ has only one real eigenvalue.  (Actually, in this case, it is easy to see that the three eigenvalues of  $\FF_{3}$ are $\{(abc)^{1/3}, (abc)^{1/3}\frac{-1 \pm i\sqrt{3}}{2}\}$.)
     
     This completes the proof of the assertions in the first paragraph of this lemma.
     
     We now consider  the assertions in the second paragraph of this lemma. Consider $H(\la) $  Suppose that $\de$ is a real root of this polynomial. Dividing $H(\la)$ by $\la-\de$ we see that the real part of the complex roots of $H(\la)$ is greater than $\de$ if and only if 
     \begin{equation}
   \phi_{1}+\phi_{2}+\phi_{3}>3\de.\label{3.38}
   \end{equation}

Suppose that $a_{1}',a_{2}', b_{1}',b_{2}',c_{1}',c_{2}'$ are all strictly positive. It follows from (\ref{3.30}) and the fact that $\phi_{1} \phi_{2} \phi_{3}\det A$ is a root of $|\Phi K-\la I|$, that $\de<\phi_{1} \phi_{2} \phi_{3}\det K$. Furthermore, since $| K|\le 1$, to prove (\ref{3.38}) it suffices to show that 
       \begin{equation}
   \phi_{1}+\phi_{2}+\phi_{3}\ge3\phi_{1} \phi_{2} \phi_{3}.\label{3.39}
   \end{equation}   
   Since $  \phi_{1}+\phi_{2}+\phi_{3}\ge3(\phi_{1} \phi_{2} \phi_{3})^{1/3}$  
  and    $\phi_{1}, \phi_{2}, \phi_{3}$ are all less than 1 we see that (\ref{3.39})
   is satisfied.  Thus, when $a_{1}',a_{2}', b_{1}',b_{2}',c_{1}',c_{2}'$ are all strictly positive and $\Phi A_-$ has only one real eigenvalue, the real part of the complex eigenvalues is greater than the real eigenvalue.   It is easy to see that this argument also works when $A_{-}$ has the form of $\FF_{1}$, $\FF_{2}$ or $\FF_{3}$.
   \qed
      
 \noindent {\bf Proof of Theorem \ref{theo-1.1} for matrices with negative off diagonal entries} This follows immediately from Lemmas \ref{general1} and \ref{lem-2.2}.\qed

 \section{$\bf3\times 3$ matrices with positive off diagonal entries}\label{sec-5}
 
  Consider the matrix

\begin{equation}
A= 
\(\begin {array}{ccc}
1&a_1&c_2\\
a_2&1&b_1\\
c_1&b_2&1
\end{array}\),\label{2.7q}
\end{equation}
in which all the entries are greater than or equal to zero.  
When $\det A>0$,
\begin{equation}
A^{-1}= \frac{1}{| A|}
\(\begin {array}{ccc}
1-b_{1}b_{2}&c_2b_{2}-a_{1}&a_{1}b_{1}-c_{2}\\
b_{1}c_{1}-a_{2}&1-c_{1}c_{2}&a_{2}c_{2} -b_1\\
a_{2}b_{2}-c_1&a_{1}c_{1}-b_2&1-a_{1}a_{2}
\end{array}\).\label{3.1bq}
\end{equation}

 \noindent {\bf Proof of Theorem \ref{theo-1.1} for matrices with positive off diagonal entries }This follows from the next lemma: 
 
    \begin{lemma}\label{theo-6.1} Let $A$, in (\ref{2.7q}),  be the kernel of a permanental vector $\th=(\th_{1},\th_{2},\th_{3})$. Then $A$ is either   diagonally equivalent to a symmetric positive definite matrix or $A^{-1}$ is an $M$ matrix.
    
 Furthermore, if   one of the off diagonal terms of $A$ is equal to 0, then $A$ is diagonally equivalent to a symmetric matrix.

 \end{lemma}

\noindent{\bf Proof }    As in (\ref{2.19}), but with $\Ga$ replaced by $A$, we have
\begin{equation}
  \Phi(\al_1, \al_2, \al_{3}) =\frac{1}{ |I+ \al A|^{ \beta}}.
   \end{equation}
  Therefore, as in (\ref{2.26})
  \begin{equation}
     \Phi_{(1,2)}(\al_{1}, \al_{2}) =\frac{1}{| I+  \al^{(2) }\Ga^{(2) } |^{\bb}}.\label{2.26a}
   \end{equation} 
   where, by (\ref{2.26-})
  \be \Ga^{(2)}= \left(
\begin{array}{cc}
1-vc_1c_2 & a_1-vc_2b_2\\
a_2-vb_1c_1 & 1-vb_1b_2
\end{array} \right)
\ee
and where $v=\displaystyle\frac{u_{3}}{1+u_{3}}$. It follows from Lemma \ref{lem-2.1} that $\Ga^{(2)}$ is the kernel of a permanental vector. Therefore, by (\ref{vj.2})
   \be (a_1-vc_2b_2)(a_2-vb_1c_1)\geq 0\label{5.6a}. \ee 
Suppose none of the off diagonal entries of $A$ are equal to 0.  The inequality in (\ref{5.6a}) holds for all $v \in (0,1)$. Therefore, either  
 \be a_1\ge c_2b_2  \quad \mbox{and}\quad a_2\ge b_1c_1\label{5.7} \ee
 or there exists a $v_{0}\in (0,1)$ such that 
\be 
a_1-v_{0}c_2b_2=a_2-v_{0}c_1b_1=0.\label{5.8}
\ee
It follows from (\ref{5.8}) that $a_{1}=v_{0}c_2b_2$ and $a_2=v_{0}c_1b_1$, or, equivalently, that
\begin{equation}
   a_{1}b_{1}c_{1}=a_{2}b_{2}c_{2}.\label{5.9q}
   \end{equation}
If  (\ref{5.9q}) holds, it follows from  Lemma \ref {lem-7.1}  that       the matrix $A$ is diagonally equivalent to a symmetric matrix.

We repeat this argument twice considering $     \Phi_{(3,2)}(\al_{1}, \al_{3})$ and $     \Phi_{(3,2)}(\al_{2}, \newline\al_{3})$.  If (\ref{5.8}) holds then we get comparable equalities when we consider $     \Phi_{(3,2)}(\al_{1}, \al_{3})$ and $     \Phi_{(3,2)}(\al_{2}, \al_{3})$. However if (\ref{5.7}) holds we also get 
\begin{equation}
   b_{1}\ge a_{2}c_{2}, \quad  b_{2}\ge a_{1}c_{1},\quad  c_{1}\ge a_{2}b_{2}\quad\mbox{and}\quad c_{2}\ge a_{1}b_{1}.\label{5.9}
   \end{equation}
    It follows from (\ref{5.7}) and (\ref{5.9}) that $A$ is an $M$-matrix if it is invertible, or equivalently, $|A|>0$.

 However, if  $| A|=0$, $A$ does not have an inverse and the consideration of whether $A^{-1}$ is an $M$-matrix is meaningless.   Therefore we must show that when  (\ref{5.7}) and (\ref{5.9}) hold   and $A$ is   not diagonally equivalent to a symmetric positive definite matrix then $|A|>0$.
 
 We need only consider the case in which $A$ is not symmetric. 
Without loss of generality we can consider that 
\begin{equation}
A= 
\(\begin {array}{ccc}
1&a &c_{1} \\
a &1&b \\
c _{2}&b &1
\end{array}\).\label{2.7jj}
\end{equation}
where $a, b, c >0$, $a^{2}=a_{1}a_{2}$, $b^{2}=b_{1}b_{2}$,   $ c _{1}c _{2}=c^{2}$ and  $c_{1}\ne c_{2}$.  
Let $d$ be such that $c_1+c_2=dc$. Obviously $d >2$.

 Since $|A|=0$
we have \be 1-(a^2+b^2+c^2)+abcd=0\label{5.13} \ee
We consider $d$ in (\ref{5.13})  as a function of $a, b, c $,  i.e.,
 \be d(a, b, c)=\frac{(a^2+b^2+c^2)-1}{abc} .\ee
Note that the gradient of $d(a,b,c)$
 \bea 
\lefteqn{\nabla  d(a,b,c)=\frac{1}{(abc)^2}\(  bc(1+a^2-b^2-c^2),\right.\label{grad}}\\
 & &\quad\qquad\qquad\nn \left.  ac(1+b^2-a^2-c^2), ab(1+c^2-a^2-b^2)\).
\eea
The inequalities in (\ref{5.7}) and (\ref{5.9}) hold when $A$ is not diagonally equivalent to a symmetric matrix. (The argument we give does not require that $|A|>0$.) When they hold we see that the components of $\nabla  d(a,b,c)$ are all greater than or equal to 0. For example, $a\ge bc$ implies that
\begin{equation}
   (1+a^2-b^2-c^2)\ge (1+(bc)^2-b^2-c^2)=(1-b^{2})(1-c^{2})\ge 0.
   \end{equation}
Note that $d(1,1,1)=2$.  Therefore, since   $a, b,c$  are all less than or equal to 1, $d(a,b,c)\le 2$. This contradiction shows that here are no permanental vectors with $|A|=0$ other than those with kernels that are diagonally equivalent to a symmetric matrix. 

 To show that if   one of the off diagonal terms of $A$ is equal to 0, then $A$ is diagonally equivalent to a symmetric matrix we consider (\ref {5.6a}).
Suppose $a_1=0$ then  either $a_2=0$ or one of $b_2, c_2$ is equal  to $0$. In these cases it follows from Lemma \ref{lem-7.1} that $A$ is effectively equivalent to a symmetric matrix.  Using  $\Phi(2,3)$ and $\Phi(1,3)$ we come to the same conclusion for all the other   ways  one or more of the off diagonal terms of $A$ can be  equal to 0. \qed

\begin{example} {\rm It seems worthwhile to point out that there are many symmetric matrices with positive entries that have determinant 0.   All $3\times 3$ symmetric matrices of the form of $D$ in (\ref{3.1}), with, $|a|\le 1$, $|b|\le 1$ and $|c|\le 1$ and with $|D|=0$  have the form
\begin{equation}
\SS_{\pm}(x,y)=\(\begin {array}{ccc}
1&\sin x &\cos y  \\
\sin x  &1&\sin (x\pm y) \\
\cos y &\sin(x\pm y) &1
\end{array}\), \label{5.16a}
\end{equation}
for any $x$ and $y$ which $\sin x,\cos y , \sin (x+ y)$ or $\sin x,\cos y , \sin (x- y)$  are greater than or equal to zero. To get this we note that $|D|=0$ implies that 
\begin{equation}
   c=ab\pm\((1-a^{2})(1-b^{2})\)^{1/2}.\label{5.17s}
   \end{equation}
If we take  $a=\sin x$, and $b=\cos y$ and solve for $c$ we get (\ref{5.16a}).

We now ask for what values of $x$ and $y$ is the adjugate, (also called the adjoint) of $\SS_{\pm}(x,y)$
a singular $M$-matrix.  (I.e. even though the matrix is not invertible, the adjugate has negative, including 0,  off diagonal elements.) Referring to $D$, and noting (\ref{3.1bq}), this is equivalent to asking for what values of $x$ and $y$ are 
\begin{equation}
   c\ge ab,\quad a\ge bc,\quad \mbox{and}\quad b\ge ac.\label{5.18}
   \end{equation}
   To achieve the first inequality in (\ref{5.18}) we must use the plus sign in (\ref{5.17s}) which gives $c=\sin(x+y)$.  This implies that $a$, which satisfies an analogue of (\ref{5.17s}), satisfies
     \begin{equation}
   a=bc-\((1-b^{2})(1-c^{2})\)^{1/2}. 
   \end{equation} 
To get the second inequality in (\ref{5.18}) we can take $c=1$ and $a=b$, which can be achieved by taking $x\in [0,\pi/2]$ and $y=(\pi/2)-x$, or $a=c$ and $b=1$ which can be achieved by taking $x\in [0,\pi/2]$ and $y=0$. In either case we get
  matrices of the form  
\begin{equation}
\AA(a)=\(\begin {array}{ccc}
1&a &a  \\
a  &1&1\\
a&1  &1
\end{array}\) \qquad a\in [0,1],\label{5.17}
\end{equation}and the matrices that can be obtained from them by interchanging their rows and columns. The adjugate of $\AA(a)$ is 
\be  \AA'(a)=  
\(\begin {array}{ccc}
0& 0&0\\
0&1-a^{ 2}&-(1-a^{ 2})\\
0&-(1-a^{ 2})&1-a^{ 2}
\end{array}\).
\end{equation}}
\end{example}
    The next lemma is an analogue of Lemma \ref{general1}   when the off diagonal elements of the kernel are all greater than or equal to zero. It also shows  that  the necessary condition  in Lemma \ref{lem-2.2}  is satisfied in this case.

     \bl \label{general2}
Let $ A_+$ be the   matrix  
\begin{equation}
 A_+= 
\(\begin {array}{ccc}
1& a'_1& c'_2\\
a'_2&1&b'_1\\
c'_1&b'_2&1
\end{array}\).\label{2.7wss}
\end{equation}
in which $a_{1}',a_{2}', b_{1}',b_{2}',c_{1}',c_{2}'$ are all greater than or equal to zero, with  $a_{1}'a_{2}',b_{1}'b_{2}',c_{1}'c_{2}'$ all less than or equal to 1. If $A_{+}$ is not   diagonally equivalent to a symmetric matrix and $ A_+^{-1}$ is an $M$-matrix,  there exists a  diagonal matrix $\Phi$, with strictly positive diagonal entries, such that $\Phi A_+$ has only one real eigenvalue. 

  Furthermore, the real part of the complex eigenvalues of $\Phi A_{+}$  is less than the real eigenvalue.  
\el

\noindent{\bf Proof of Lemma   \ref{general2} } Assume, to begin, that  $a_{1}',a_{2}', b_{1}',b_{2}',c_{1}',c_{2}'$ are all strictly positive.
 We use the notation of the proof of Lemma \ref{general} except that we consider the associated symmetric matrix
  $K_{+}$ be the matrix
 \begin{equation}
 K_{+}= 
\(\begin {array}{ccc}
1&a&c\\
a&1&b\\
c&b&1
\end{array}\), \label{3.1b}
\end{equation}
where $a_{1} 'a_{2}'=a^{2}$, $ b_{1}' b_{2}'=b^{2}$ and $ c_{1} 'c_{2}'=c^{2}$. 
  It follows from Remark \ref{rem-2.2} that $\det K_{+}>0$.  As in the proof of   Lemma \ref{general} 
 \bea
\lefteqn{ \det\(\Phi K_{+}-\la I\)}\label{7.20sos}\\
&=&   -\la^3+\la^2\( \phi_1+ \phi_2+ \phi_3 \)\nn  \\  
& &\hspace{.2in} -\la \(  \phi_1 \phi_2(1-a^2)+ \phi_2 \phi_3(1-b^2)+ \phi_1 \phi_3(1-c^2)\)+ \phi_1 \phi_2 \phi_3 \det K\nn\\
 &=&(\la - 1)^2(\phi_1\phi_2\phi_3 \det K_{+}-\la ).\nn
      \eea
      Note that 
      \begin{equation}
   \frac{d^{2}}{d\la^{2}} \det\(\Phi K_{+}-\la I\)=-6\la+2(   \phi_{1}+\phi_{2}+\phi_{3}).\label{4.45}
   \end{equation}
 Using the fact that \begin{equation}
   \phi_{1}+\phi_{2}+\phi_{3}>3
   \end{equation}we see  that $\det\(\Phi K_{+}-\la I\)$ has a local minimum of 0 at $\la=1$.
   
 Consider 
 \be
 \det\(\Phi  A_{_+}-\la I\)
 \ee
By Lemma  \ref{lem-2.3a} it suffices to take $ A_{+}$ to be 
\begin{equation}
\(\begin {array}{ccc}
1& a & c \\
 a &1& b_1\\
 c & b _2&1
\end{array}\).\label{2.7hh2}
\end{equation}
where   $b _{1}b _{2}=b ^{ 2}$. We have
\be 
 \det K_{+}=1-a^{2}-b^{2}-c^{2}+2abc 
 \ee 
 and
  \be 
 \det  A _{+}=1-a^{2}-b^{2}-c^{2}+ac (b_{1}+b_{2}).
 \ee 
 Since $b _{1}b _{2}=b ^{ 2}$ we see from the first equality in (\ref{7.20sos}) that    
 \begin{equation}
  \det\(\Phi  A_{ +}-\la I\)-\det\(\Phi K_{+}-\la I\) =ac (b_{1}+b_{2}-2b).\label{3.50}
   \end{equation}  
 Since $\det\(\Phi K_{+}-\la I\)$ has a local minimum of 0 at $\la=1$,  when $b_{1}\ne b_{2}$, $\det\(\Phi  A _{+} -\la I\)$ has a local minimum that is strictly positive at $\la=1$. This implies that when $b_{1}\ne b_{2}$, $\det\(\Phi  A_{+}-\la I\)$ has only one real root. This is equivalent to the statement of this lemma. 
 
 Suppose now that one of the terms  $a_{1}',a_{2}', b_{1}',b_{2}',c_{1}',c_{2}'$ is equal to zero. Let's say it is $a_{1}'$. Then since
\begin{equation}
A_{+}^{-1}= \frac{1}{\det A_{+}}
\(\begin {array}{ccc}
1-b'_{1}b'_{2}&c'_2b'_{2}-a'_{1}&a'_{1}b'_{1}-c'_{2}\\
b'_{1}c'_{1}-a'_{2}&1-c'_{1}c'_{2}&a'_{2}c'_{2} -b'_1\\
a'_{2}b'_{2}-c'_1&a'_{1}c'_{1}-b'_2&1-a'_{1}a'_{2}
\end{array}\),\label{3.1b2}
\end{equation}
 we see that either $b_{2}'$ or $c_{2}'$ must be zero. Suppose it is $b_{2}'$. Then $A_{+}$ in (\ref{2.7hh2}) has the form
 \begin{equation}\wt A_{+}=
\(\begin {array}{ccc}
1& 0& c'_2\\
a'_2&1&b'_1\\
c'_1&0&1
\end{array}\).\label{2.7io}
\end{equation}
  and therefore it is equivalent to the symmetric matrix 
 \begin{equation}\ov A_{+}=
\(\begin {array}{ccc}
1& 0& c \\
0&1&0\\
c &0&1
\end{array}\).\label{2.7iof}
\end{equation}

We now show that when $\Phi A_{+}$ has only one real eigenvalue, the real part of the complex eigenvalues is less than the real eigenvalue. Let $\de$ denote the real eigenvalue.  Referring to (\ref{3.38}) we see that we must show that
   \begin{equation}
   \phi_{1}+\phi_{2}+\phi_{3}<3\de.\label{3.38q}
   \end{equation}
  It follows from (\ref{3.50}) that the roots of  $ \det\(\Phi  A_{ +}-\la I\)$ are greater than the roots of $\det\(\Phi K_{+}-\la I\)$. Therefore, by (\ref{7.20}),   $\de>   \phi_{1} \phi_{2} \phi_{3}\det K_{+}$. Thus to obtain (\ref{3.38q}) it suffices to show that 
 \begin{equation}
      \phi_{1}+\phi_{2}+\phi_{3}\le 3\phi_{1} \phi_{2} \phi_{3}\det K_{+}.\label{3.38q2}
   \end{equation}
By (\ref{4.45})    the second derivative of $ |\Phi K_{+}-\la I|$ is  negative when $\la>(\phi_{1}+\phi_{2}+\phi_{3})/3$.   Considering the graph of  $ |\Phi K_{+}-\la I|$ we see that this implies that   the single real root of $ |\Phi K_{+}-\la I|$ is greater than $(\phi_{1}+\phi_{2}+\phi_{3})/3$; hence we have (\ref{3.38q}).\qed 

 \section{ Premanenetal vectors with pairwise independent components.}\label{sec-6}

   \noindent {\bf Proof of Theorem \ref{theo-1.2} }
  We first consider the case in which
 $\th\in R_{+}^{3}$.     
    Let $\Ga$ be the kernel of $\th$.   It is enough to consider the case in which all the diagonal elements of $\Ga$ are equal to one. Since $\th$ has  pairwise independent components we know that (\ref{4.2v}) holds. In this case either $\det \Ga=\Ga(1,1)\Ga(2,2)  \Ga(3,3)=1$ or else 
\begin{equation}
\Ga = 
\(\begin {array}{ccc}
1&0&a\\
b&1&0\\
0&c&1
\end{array}\)\label{3.1w}
\end{equation}
with $abc \neq 0$, or $\Ga^{T}$ is equal to this matrix. (It is obvious $\Ga$ must contain three zeros off the diagonal. Any configuration other than (\ref{3.1w})  or its transpose has determinant equal to 1.) 

Suppose the off diagonal elements of $\Ga$ are positive. It is obvious from (\ref{3.1bq}) that $\Ga^{-1}$, if it exists, is not an $M$-matrix. Therefore, by Theorem \ref{theo-1.1}, $\Ga$ is diagonally equivalent to a symmetric matrix, which implies that $abc=0$. If the off diagonal elements of $\Ga$ are negative Theorem \ref{theo-1.1} again implies that $\Ga$ is diagonally equivalent to a symmetric matrix, which again implies that $abc=0$.  Therefore Theorem \ref{theo-1.2} holds for $\th\in R^{3}_{+}$.

We next consider a generalization of (\ref{3.1w}) to $n \times n$ matrices,  $\wt\Ga_{n}$, $n \geq 3$. These are matrices for which
 \begin{equation}
   |I+\al\wt  \Ga_{n}|=\prod_{i=1}^{n}(1+\al_{i})+\CC(\wt  \Ga_{n})\prod_{i=1}^{n} \al_{i} ,\label{3.6}
   \end{equation}  
in which $\CC(\wt  \Ga_{n})$ is a real valued function of the components of $\wt  \Ga_{n}$. We have the following lemma:

\bl \label{lem-3.2}
 Let $\th\in R_{+}^{n}$,  $n \geq 3$, be a $\bb$-permanental vector, with pairwise independent components and kernel $\wt  \Ga_{n}$.  Then   $\CC(\wt  \Ga_{n})=0$, in which case the components of $\th$ are independent.
\el

\Proof     We show in the beginning of this  section that this lemma is true when $n=3$. For $n>3$ the Laplace transform of $\th$ is  
\be
 \Phi(\al_{1},\ldots,\al_{n})=\left|\prod_{i=1}^{n}(1+\al_{i})+\CC(\wt  \Ga_{n})\prod_{i=1}^{n} \al_{i}\right|^{-\bb}.
\ee
 
Let $\al_i =1$ for all $4 \leq i \leq n$. By Lemma \ref{lem-2.2}  
 \bea
 \wt \Phi (\al_{1},\al_{2},\al_{3})& :=& \frac{\Phi(\al_{1},\al_{2},\al_{3},1,\ldots,1)}{ 
 \Phi(0 ,0,0,1,\ldots,1) } \label{5.4}\\
  &=&\nn\left|(1+\al_1)(1+\al_2)(1+\al_{3})+ { \CC(\wt  \Ga_{n})\over 2^{n-3}}\al_1\al_2\al_{3}  \right|^{-\bb} 
   \eea
   is a Laplace transform of a random variable in $R^{3}_{+}$. Furthermore, the form of the right-hand side of (\ref{5.4})  shows that   $(\th_{1},\th_{2},\th_{3})$
   is a permanental vector with a kernel of the form of (\ref{3.1w}), (or its   transpose).
   We show in the beginning of this section that for such a vector we must have    $ \CC(\wt  \Ga_{n})=0$.\qed

 	  \noindent{\bf Proof of Theorem \ref{theo-1.2} continued }
 The proof is by induction.    We show in the beginning of this section that Theorem \ref{theo-1.2} holds for  for $n=3$. Let $n\ge 3$ and assume that the theorem holds for all  $m <n$. Let $\th$ be a $\bb$-permanental process in $R^{n}_{+}$ with pairwise independent components and kernel $G$. This theorem follows from Lemma \ref{lem-3.2} once we show that the determinant $|I+\al G|$  for this process has the form of (\ref{3.6}).
 
 Suppose it does not. Then the it must contain a term of the form
 \begin{equation}
   \CC'(G)\prod_{j=1}^{k}(1+ \al_{i_{j}}G(i_{j},i_{j}))
   \end{equation} 
 where $1\le k<n$ and $(  {i_{1}} ,\ldots, {i_{k}})$ is a proper subset of 
 $( { {1}} ,\ldots, {n})$ and $   \CC'(G)\ne 0$. In fact we know what $   \CC'(G)$ is. Let $( \al_{j_{1}} ,\ldots,\al_{j_{n-k}})$ be the elements of  $( \al_{ {1}} ,\ldots,\al_{n})$ that are not in  $( \al_{i_{1}} ,\ldots,\al_{i_{k}})$ and $G'$ be the $(n-k)\times( n-k)$ matrix obtained by removing the $   i_{1}$-th,$    \ldots,$$ i_{k}  $-th row and column from $G$. It is easy to see that
 \begin{equation}
    \CC'(G)=|I+\ga G'|
   \end{equation}
 where $I$ is the $(n-k)\times( n-k)$ identity matrix and $\ga$ is the diagonal matrix with $\(\ga_{i,i}=\al_{j_{i}}\)$,   $1\le i\le n-k$. 
 
 By the hypothesis of this theorem $|I+\ga G'|^{-\bb}$ is the Laplace transform of a permanental process in $R^{n-k}_{+}$ with pairwise independent components. By the induction hypotheses these components are independent. Therefore
 \begin{equation}
    \CC'(G)=\prod_{i=1}^{n-k}(1+\al_{j_{i}}G(j_{i},j_{i})).
   \end{equation}
 This shows that  the determinant $|I+\al G|$  for this process has the form of (\ref{3.6}). 
 \qed

 	 \section{Proof of Corollary  \ref{cor-1.1} }\label{sec-7}
 It is easy to see what we must show. For all $x,y,z\in T$ \begin{equation}
   d(x,y)\le d(x,z)+d(y,z).\label{7.1}
   \end{equation} 
   This follows if we can show that the determinant of 
   \begin{equation}
   \widehat\Ga=
\(\begin {array}{ccc}
\Ga(x,x)&(\Ga(x,y)\Ga(y,x))^{1/2}&(\Ga(x,z)\Ga(z,x))^{1/2}\\
(\Ga(x,y)\Ga(y,x))^{1/2}&\Ga(y,y)&(\Ga(y,z)\Ga(z,y))^{1/2}\\
(\Ga(x,z)\Ga(z,x))^{1/2}&(\Ga(y,z)\Ga(z,y))^{1/2}&\Ga(z,z)
\end{array}\)\label{10.1}
\end{equation}
   is greater than or equal to zero. This follows because, by (\ref{3.2k}), the $2\times 2$ principal minor of $\wh\Ga$ is      greater than or equal to zero. Therefore,  if $|\wh\Ga|\ge 0$, $\wh\Ga$ is positive definite and hence the covariance of a Gaussian vector in $R^{3}$. Considering (\ref{1.7}) we get (\ref{7.1}).
  
  By Lemma \ref{lem-1.1} the kernel of the permanental vector 
 $(P(x),P(y),P(z))$ is 
 \begin{equation}
 \Ga=
\(\begin {array}{ccc}
\Ga(x,x)& \Ga(x,y) & \Ga(x,z)  \\
 \Ga(y,x) &\Ga(y,y)&\Ga(y,z)\\
 \Ga(z,x) & \Ga (z,y) &\Ga(z,z)
\end{array}\)\label{10.1a}.
\end{equation}

Since $\Ga$ is the  kernel of a permanental vector we have $|\Ga|\ge 0$. We must show that  
\be |\Ga|\ge 0\quad\mbox{ implies that}\quad  |\wh \Ga|\ge 0.
\ee
 One of the idiosyncrasies that we must take into account is that the off diagonal elements of $\wh \Ga$ are always greater than or equal to zero whereas the off diagonal elements of $  \Ga$ may be negative.

 Suppose  $\Ga\ge 0$   and $|\Ga|\ge0$. By Theorem  \ref{theo-1.1} and Lemma \ref{lem-7.1} either $\Ga^{-1}$ is an $M$-matrix or else $\Ga$ is diagonally equivalent to a symmetric positive definite matrix. In the first case, by Lemma \ref{lem-1.1} and Remark \ref{rem-2.2}, $|\wh\Ga|>0$.
In the second case, it is easy to see from (\ref{2.14}) that $|\Ga|=|\wh\Ga|$ so we also have $|\wh\Ga|>0$.

When $\Ga$ has negative off diagonal terms   we consider two cases. The first  is that $\Ga$ is diagonally equivalent to a matrix with positive off diagonal terms, say $\Ga'$, i.e., $\Ga=D\Ga' D^{-1}$. In this case  $|\Ga|=|D\Ga' D^{-1}|$ and the matrix in (\ref{10.1}) is the same for $\Ga$ and $\Ga'$. Therefore, the argument in the previous paragraph   shows that $|\wh \Ga|>0$.

Now, suppose  $\Ga$ has negative off diagonal terms and it is not  diagonally equivalent to a matrix with positive off diagonal terms. Relabel the matrix $\Ga'$. Without loss of generality we can assume
  \begin{equation}
 \Ga'=
\(\begin {array}{ccc}
\Ga(x,x)& -\Ga(x,y) & -\Ga(x,z)  \\
 -\Ga(y,x) &\Ga(y,y)&-\Ga(y,z)\\
- \Ga(z,x) & -\Ga (z,y) &\Ga(z,z)
\end{array}\)\label{10.1b}.
\end{equation}
  By Lemma \ref{general} $\Ga'$ is diagonally equivalent to a symmetric matrix. 
 By Lemma \ref{lem-7.1}, $|\Ga'|$ is equal to     \begin{equation} 
\left|\begin {array}{ccc}
\Ga(x,x)&-(\Ga(x,y)\Ga(y,x))^{1/2}&-(\Ga(x,z)\Ga(z,x))^{1/2}\\
-(\Ga(x,y)\Ga(y,x))^{1/2}&\Ga(y,y)&-(\Ga(y,z)\Ga(z,y))^{1/2}\\
-(\Ga(x,z)\Ga(z,x))^{1/2}&-(\Ga(y,z)\Ga(z,y))^{1/2}&\Ga(z,z)
\end{array}\right|\label{7.5}
\end{equation}
 Therefore, since  $\Ga'$ is the determinant of a permanental vector $ |\Ga'|\ge 0$.    This implies  that $|\wh \Ga|\ge0$ in this case also. \qed

\section{$M$-matrices and infinite divisibilty}\label{sec-8}

  Critical in the proof of Theorem \ref{theo-1.1} is the fact  that if a matrix $\Ga$ is invertible and  $\Ga^{-1}$ is diagonally equivalent to an $M$-matrix, then $\Ga$ is the kernel of a permanental vector for all $\bb>0$. This implies that the vector is infinitely divisible. Necessary and sufficient conditions for a  vector of Gaussian squares to be infinitely divisible are due to R. Bapat and R. C. Griffiths \cite{Bapat, Griffiths}. An exposition of this material is given in \cite[Chapter 13]{book}. Eisenbaum and Kaspi, \cite[Lemma 4.2]{EK} recognize that   the sufficient condition in the Bapat--Griffiths criteria for infinite divisibility in the case of symmetric kernels also works for non-symmetric kernels. Their proof involves probabilistic considerations. Since Theorem \ref{theo-1.1} is only for $3\times 3$ matrices it seems appropriate to give a proof for finite matrices involving only linear algebra.

  As one might expect the proof mimics the proof in the symmetric case which is given in \cite[Lemma 14.9.4 and Theorem 13.2.1]{book}. In fact the more subtle revisions need to be made in the lemma. We do this next. What we prove is given
in \cite[Theorem 2.3]{B-P}. The proof in \cite{B-P}   omits all details. 

\bl \label{mequiv}
Let  $B=\{B_{i,j}\}_{1\le i,j\le p}$ be an    invertible matrix such that $B_{i,i} >0$ for all $i$ and  $B_{i,j} \leq 0$ for all $i \neq j$.
The following are equivalent:
\begin{enumerate}
\item  $B$ is an $M$-matrix; i.e., $B^{-1} \geq 0$.
\item  There exists a positive diagonal matrix $D$ such that all the eigenvalues of  $BD$ have strictly positive real parts.
\item  There exists a matrix $C\geq 0$ such that $ BD+C= \la I$ where $\la$ is greater than the spectral radius of $C$. (I.e.,  $\la >\rho $, where $\rho:=\rho(C)$ is the magnitude of the maximal eigenvalue of $C$.)
\end{enumerate}
\el

\Proof    We first show that  {\it 1.}  implies {\it 2. }
Let $D$ be a diagonal matrix with diagonal entries $D_{i,i}=\sum_{j=1}^p
\{B^{-1}\}_{i,j}$, $i=1,\ldots,n$. Since $B^{-1}\ge 0$, $BD$ has positive diagonal elements and negative off--diagonal elements. Note that \bea
  \sum_{k }  \{BD\}_{i,k}&=&\sum_k B_{i,k}D_{k,k}=\sum_k B_{i,k}\sum_j \{B^{-1}\}_{k,j}\\
&=&\sum_{j,k} B_{i,k}\{B^{-1}\}_{k,j}=\sum_j \de_{i,j}=1\nn.\eea
This shows that  that $BD$ is strictly diagonally dominant, i.e., for all $i$,
\be
\sum _{k:k\neq i}  |BD_{i,k}| <BD_{i,i}.\label{8.2}
\ee

Consider the symmetric matrix   $A=BD+(BD)^t$, where the superscript $t$ represent conjugate transpose.
To see that $A$  is positive definite, let  $u=(u_{1},\ldots,u_{p})$ be any real vector, and let  $c_{i,j}:= \{BD\}_{i,j}$.  Using  (\ref{8.2}) we have
\bea
(BDu,u)&= &  \sum_i u_i^2 c_{i,i}+\sum_{i,j:i\neq j} u_iu_j c_{i,j}  \\&> &
\sum_i u_i^2 \sum_{j:j\neq i} |c_{i,j}| - \sum_{i,j:i\neq j} |u_iu_j| |c_{i,j}| \nn \\&=& \sum_{i,j:i<j}\(u_{i}^{2} |c_{i,j}|+u_{j}^{2} |c_{j,i}|-|u_{i}||u_{j}|(|c_{i,j}|+|c_{j,i}|)\)  \nn
\eea
It follows from this that
\bea
(Au,u)&>& \sum_{i,j:i<j}\(u_{i}^{2}  +u_{j}^{2} -2|u_{i}||u_{j}|\)\(|c_{i,j}|+|c_{j,i}|\)\nn \\ &=& \sum_{i,j:i<j} ( |u_i|-|u_j|)^2\(|c_{i,j}|+|c_{j,i}|\)\geq 0.\nn \\
\eea
Since $A$ is a
  symmetric and strictly positive definite matrix its  eigenvalues are all greater than 0.

  In general, if
  $\al$ is an eigenvalue of $BD$, $\bar{\al}$ is an eigenvalue of  $(BD)^t$.  Let $\La$ be the eigenvalue matrix of $BD$, and write $BD= E \La E^t$, where $E$ is the matrix of eigenvectors of $BD$.   Then
\be A=  E \La E^t+  E \La^{t} E^t =E(\La+\La^{t})E^t\label{8.5}.
\ee
Therefore, since the eigenvalues of $A$ are all greater than 0, the eigenvalues of $BD$ have strictly positive real parts.

To show that 2. implies 3. we first
note  that  because the off diagonal elements of $BD$ are less than or equal to zero, it is easy to see that we can find a positive number $\la$ and matrix $C'\ge 0$ such that
\begin{equation}
  BD=\la I-C'.\label{8.6}
  \end{equation}
It remains to show that we can do this with   $\la>\rho(C')$.

Clearly, for any $p$-dimensional complex valued unit vector $y$,  $|C'y|\leq \rho $. Let $\wt x$ be an  eigenvector of $C'$, of unit length, such that  $C' \wt x= \ga \wt x$, where $|\ga|=\rho$.  Let $x$ be the vector with  $Re(x_i)=|Re(\wt x_i)|$  and $Im(x_i)=Im(\wt x_i) $ for all $i$.  Trivially, $|x|=|\wt x|$. Then, since $C'$ is positive, unless $x=\wt x$,  $| C'x| > |C'\wt x |=\rho $. Here we have assumed that a least one real component of $\wt x$ is not equal to 0. If that is not the case then $C'(i\wt x)=\ga (i\wt x)$ and since $i\wt x$ has real components we can apply the argument above to see that we must have $x=i\wt x$.
Consequently, there exists an eigenvector $x$ of $C'$ such that $Re(x_i) \geq 0$ for all $1\le i\le p$, with strict inequality for at least one $i$, satisfying \begin{equation}
  C'x= \ga x.\label{8.5a}
  \end{equation}

  Let $A$ be the strictly positive definite symmetric matrix defined above in this proof.  The argument used to obtain (\ref{8.5}) shows that there exists a $\la$ such that
  \begin{equation}
  2\la> \ga+\bar\ga .\label{8.8}
  \end{equation}
    Let $ \ga=a+ib$, in which $a$ and $b$ are real.  By (\ref{8.8}) there exists an $\ep>0$ such that  $\la = a+ e$.

    For any $\de>0$, $BD+(C'+\de I)=(\la +\de)I$. The spectral radius of $C'+\de I$ is   $|(a+\de)^2+b^2|^{1/2}$. It is easy to see that for sufficiently large $\de$, say $\de_{0}$, this is smaller than $(\la +\de_{0})^2$. Taking $C=(C'+\de_{0} I)$ we get 3.

We conclude the proof by showing that 3. implies 1.
Clearly the absolute values  of  all the  eigenvalues of $C/\la$  are less than one.  Therefore we can write
\be (BD)^{-1}=\frac{1}{\la}\left(I- \frac{C}{\la}\right )^{-1}=\frac{1}{\la} \sum _{k=0}^\infty \left(\frac{C}{\la}\right)^k
\ee
to see that $(BD)^{-1}=D^{-1}B^{-1}\ge 0$. Since $D$ is a positive diagonal matrix, 1. follows by multiplying $D^{-1}B^{-1}$ by $D$.
\qed

In the next theorem, rather than focus on infinite divisibility, we  characterize the  kernels that give $\bb$-permanental processes for all $\bb>0$.

\bt \label{mmtrx}
Let $B$ be an $n \times n$ invertible matrix and let $A= B^{-1}$:
  If there exists a signature matrix $\NN$ such that $\NN B\NN$
  is  an $M$-matrix, then  $  \NN A \NN  $ is a kernel of a  $\bb$-permanental vector for all $\beta >0$. (I.e., $|I+\al      \NN A \NN  |^{-\beta}$ is the Laplace transform of a random positive  vector for all $\beta >0$.)
\et

\Proof  One of the main steps in this proof  is that for any matrix $U$ with spectral radius less than $1$,
\bea
\log(|I-U|)=  \sum_{n=1}^\ff \frac{  \mbox{\rm trace}\{U^n\}}{n}.\label{8.10}
\eea
This is explained on \cite[page 562]{book}, for matrices that are closely related to symmetric matrices. However the proof goes through exactly as written  for any matrix  with the spectral radius less than $1$.

We proceed with the proof.   To simplify the  notation  we might as well assume that $B$ itself is a $M$-matrix, i.e. that we can take $\NN=I$. Then if we can show that there exist two strictly positive diagonal matrices  $D_{1}'$ and $D_{2}'$ such that $D'_{1}AD'_{2}$ is a kernel of a  $\bb$-permanental vector for all $\beta >0$,  the conclusion of 1. follows from Lemma \ref{lem-2.3}.

Since $B$ is an $M$-matrix, it follows from Lemma \ref{mequiv}, 3. that  for some positive diagonal matrix $D$, $\wt B =DB=\la I -C$ for some positive matrix $C$ and  $\la$ with the property  that $\la>\rho(\wt B) $.
Let  $\wt A  = (\wt B)^{-1}$ and $S$ be a diagonal matrix with diagonal entries $0\le s_{i}\le 1$, $i=1,\ldots,n$. Then for any diagonal matrix $\al$ with diagonal entries $\al_{i}$, $i=1,\ldots,n$, we can can find a positive number $t$ and a diagonal matrix $S$ such that $\al=t(I-S)$. We have
\bea
|I+\al \wt A|&=& |I+t(I-S)\wt A|  \\&=&|\wt A| |\wt B+ t(I-S)\ |\nn\\&=& |\wt A| |(\la +t)I-(C+tS)|\nn\\ &=& |\wt A|(\la+t)\left|I-\frac{(C+tS)}{(\la +t) }\right |\nn
\eea

Since $S$ has all entries less than $1$ and all eigenvalues of $C$ are less than $\la$ in magnitude, it follows from (\ref{8.10}) that
\bea
\log(|I+\al \wt A|)=  \log (\wt A)+ \log(t+\la) + \sum_{k=1}^\ff \frac{\mbox{\rm trace}\{(C+tS)^k\}}{k(\la+t)^k} \label{8.12}
\eea
for sufficiently large $\la>0$.

Since $C$ is positive it is obvious that the expansion of this power series in  about $S=0$ has non-negative coefficients. By \cite[Lemma 13.2.2]{book},  this implies that $|I+\al \wt A|^{-\beta}$ is the Laplace transform of a positive random vector for all $\beta$. Since $\wt A=D^{-1}B$, the comments in the second paragraph completes the proof. \qed 

\begin{remark}{\rm \label{VJ_counterex}   \cite[Proposition 4.6]{VJ}   states that a necessary condition for a kernel $\Ga$ to   define a permanental vector in (\ref{1.8}), for all $\bb>0$,  is that $\Ga$ and all matrices obtained from $\Ga$ by multiplying its rows by non-negative real numbers numbers, have only real non-negative eigenvalues. This is not correct even for symmetric matrices since kernels of Gaussian squares that are in class 1. and not in class 2. satisfy this condition but,  since they are not in class 2., they are not kernels of permanental vectors for all $\bb>0$. In fact, considering Lemmas  \ref{general1}  and \ref{general2}, a correct statement is:    If a  $3\times 3$ matrix $\Ga$ is the kernel of a permanental vector, and it, and all matrices obtained from it by multiplying its rows by non-negative real numbers   have only real non-negative eigenvalues,  then $\Ga$ is diagonally equivalent to a symmetric positive definite matrix. }
\end{remark}

    \def\wh{\widehat}
\def\ol{\overline}

 \def\noopsort#1{} \def\printfirst#1#2{#1}
\def\singleletter#1{#1}
            \def\switchargs#1#2{#2#1}
\def\bibsameauth{\leavevmode\vrule height .1ex
            depth 0pt width 2.3em\relax\,}
\makeatletter
\renewcommand{\@biblabel}[1]{\hfill#1.}\makeatother
\newcommand{\bysame}{\leavevmode\hbox to3em{\hrulefill}\,}

\end{document}